\documentclass[12pt]{article}
\usepackage{latexsym,amssymb,amsmath,amsfonts,latexsym,amscd,amsthm}

\newtheorem{Thm}{Theorem}

\newtheorem{Lem}{Lemma}

\newtheorem{Prop}{Proposition}
\newtheorem{Def}{Definition}

\newtheorem{Rem}{Remark}
\begin{document}
                               
\title{\bf\ Invariant tubular neighborhood theorem for affine varieties.
\footnote{Research supported by MURST 60\%}}
\author{ M. Boratynski\\
\date{}
\small Dipartimento di Matematica\\
\small via Orabona 4\\
\small 70125 Bari, Italy\\
\small boratyn@dm.uniba.it}
\maketitle
\begin{abstract}
The  aim of this note is to prove the
algebraic geometry analogue of the
Invariant tubular neighborhood theorem which concerns the actions of
compact Lie groups on smooth manifolds.
\end{abstract}
\section*{Introduction.}
Let $G$ be a compact Lie group acting on a smooth manifold $M$. Suppose
$A\subset M$ is a smooth invariant closed manifold.
 Then
 it is well known (\cite{B:4})
 that $A$ admits an open invariant tubular neighborhood
 in $M$. This means that there exist a (smooth) $G$-vector bundle $E$
 on $A$ and an equivariant diffeomorphism $\phi\colon E \to M$ onto some
 open neighborhood of $A$ in $M$ such that the restriction of $\phi$ to the
 zero section of $E$ is the inclusion of $A$ in $M$.
 In this paper we prove in the case of linear actions of reductive groups
 on affine varieties the   following analogue of
 the Invariant tubular neighborhood theorem.\\ 
 \textit{Theorem:}
 Let $G$ be a reductive group acting linearly on an affine
 space $A^{n}_{k}$  (ch $k$=0). Suppose $X\subset
 A^{n}_{k}$ is an affine smooth $G$-invariant subvariety. Then there exist
 (an algebraic) $G$-vector bundle $E$ on $X$, a $G$-invariant open
 $U\subset E$ which
 contains the zero section of $E$ and $\phi\colon U\to A^{n}_{k}$  such that
 $\phi $ is etale, equivariant and  its restriction to the zero section of
 $E$ is the inclusion of $X$ in $A^{n}_{k}$.
  The well known Luna Slice Theorem(\cite{L:1}) concerns the case when $X$ is
  the closed orbit and states that   $\phi$ can be chosen to be
   strongly etale.This does not seem to be true in the general case.
\smallskip
\par  In the sequel all the varieties considered are over an algebraically
closed field $k$ with ch $k$=0. $G$ will always denote a reductive
 group over k.
\section*{Preliminaries.}
\indent We shall
 need some elementary results about rational modules and modules
of differentials in the equivariant context.
\par Let $R$ be an affine $k$-algebra with a rational action of $G$ and
let $M$ be an $R$-module such that $G$ acts on $M$ rationally. 
\begin{Def}
$M$ is called an $R-G$
module if $g(rm)=g(r)g(m)$ for all $g\in G,r\in R$
and $m\in M$.
\end{Def}
For the proof of the following Proposition we refer to (\cite{B:1}).
\begin{Prop}\label{P:1}       Let
\begin{equation}\label{E:a}
0\to M^{'}\to M\to M^{''}\to 0
\end{equation}
be an exact sequence of $R-G$ 
modules with $M^{''}$ finitely generated. If (\ref{E:a}) splits 
as a sequence of $R$-modules then it splits as a sequence of $R-G$ modules.
\end{Prop}
\par Let $\Omega_{R/k}$ denote the module of differentials. One can check
easily that the formula $g(dr)=d(g(r))$ for $g\in G,r\in R$ defines 
an  $R-G$ module structure on  $\Omega_{R/k}$.
\par Suppose now that $G$ acts linearly on $k[X_{1},X_{2},\ldots X_{n}]$
and $I$ is a $G$-invariant ideal of $A=k[X_{1},X_{2},\ldots X_{n}]$.
Then the induced action of $G$ on $R=A/I$ is rational and $R$  is
obviously an affine $k$-algebra.
\begin{Prop}\label{P:2}
Let $R$ be smooth over $k$.Then the first fundamental exact sequence 
(\cite{H:1})
\begin{equation}\label{E:b}
0\to I/I^{2}\to \Omega_{A/k}\otimes A/I\to \Omega_{R/k}\to 0
\end{equation}
splits as a sequence of $R-G$ modules.
\end{Prop} 
\begin{proof}
Note that all the maps in (\ref{E:b}) are $G$-invariant.
Moreover (\ref{E:b}) is a splitting sequence of $R$-modules
since $R$ is smooth.
We can now apply\linebreak Proposition \ref{P:1}.
\end{proof}

\smallskip\par
In what follows $I$ is an ideal of a ring $A$ (commutative and noetherian)
and we put $\widehat{A}=\underleftarrow{lim}A/I^{n}$

\begin{Lem}(\cite{B:3})
Let $I\subset A$ be an ideal and let $\mathfrak{m}\supset {I}$ be a maximal
ideal of $A$. Put $\widehat{\mathfrak{m}}=\mathfrak{m}\widehat{A}$.
Then the map $A_{\mathfrak{m}}\to \hat{A}_{\widehat{\mathfrak{m}}}$
induces the isomorphism of the completion of $A_{\mathfrak{m}}$ with respect
to the $\mathfrak{m}A_{\mathfrak{m}}$-adic topology with the completion of 
$\hat{A}_{\widehat{\mathfrak{m}}}$ with respect to the
$\widehat{\mathfrak{m}} \hat{A}_{\widehat{\mathfrak{m}}}$-topology. 
\end{Lem}
We immediately obtain the following
\begin{Prop}\label{P:3}
Let $G$ act rationally on affine $k$-algebras $A$ and $B$ with $G$-invariant
ideals {I} and {J} respectively.
Suppose  $d\colon A\to B$ is a $G$-invariant homomorphism which for all
$n$ induces an isomorphism $A/I^{n}\simeq B/J^{n}$.Then the set of all
the points of $Spec\, B$ at which the induced
map $Spec\, B \to Spec\, A$ is etale is 
an open, $G$-invariant subset containing                                                 
all the closed points $[\mathfrak{m}]$ such that
$\mathfrak{m}\supset {J}$
\end{Prop}
 \begin{proof}
It suffices to note that the morphism $Spec\, B\to Spec\, A$ is etale at
 $x$=[$\mathfrak{m}$] if and only if $f$ induces the isomorphism of the
 completions of the corresponding local rings (\cite {M:1}).
 Moreover the "etale" property is local and $G$-invariant for formal
 reasons since $d$ is $G$-invariant.
 \end{proof}
 \section*{Main Theorem.}
 \begin{Thm}Let $G$ 
 act linearly on an affine
 space $A^{n}_{k}$. Suppose $X\subset
 A^{n}_{k}$ is a closed, affine and smooth $G$-invariant subvariety. Then
  there exist
 (an algebraic) $G$-vector bundle $E$ on $X$, a $G$-invariant open
 $U\subset E$ which
 contains the zero section of $E$ and $\phi\colon U\to A^{n}_{k}$  such that
 $\phi $ is etale, $G$-invariant and  its restriction to the zero section
 of $E$ is the inclusion of $X$ in $A^{n}_{k}$.
 \end{Thm}
 \begin{proof}
 We put $A=k[X_{1},X_{2},\ldots X_{n}]$ and $R=A/I$ where $I$ is a $G$-
 invariant ideal of $A$ which corresponds to $X$.
 Then by Proposition \ref{P:2} there exists a
 $R-G$ homomorphism $\gamma \colon \Omega_{A/k}\otimes A/I \to I/I^{2}$ 
 which splits $I/I^{2}\to \Omega_{A/k}\otimes~A/I$.
 For $f\in A$ we put $d_{1}(f)=\gamma (\overline{df} )$ where 
 $\overline{df}$ denotes the image of $df$ in $ \Omega_{A/k}\otimes A/I$.
 Thus  $d_{1}\colon A\to  I/I^{2}$ is $G$-invariant and its restriction
 to $ I/I^{2}$ coincides with the natural homomorphism 
 $I\to I/I^{2}$. We denote by $d_{0}$ the natural homomorphism $A\to A/I$.
 \par The association $X_{i}\to d_{0}(X_{i})$+ $d_{1}(X_{i})$ 
 for $i=1,2,   \dots ,n$ defines \linebreak 
 $d\colon A \to S(I/I^{2})$ where $S(I/I^{2})$ denotes the
 symmetric
 algebra of the $R$-module $I/I^{2}$ with the obviously induced action of $G$.
 It turns out that
 $d\colon A \to S(I/I^{2})$
 is a $G$-invariant $k$-algebra homomorphism.\par Let
 $d_{i} \colon A \to S^{i}(I/I^{2})=I^{i}/I^{i+1}$ denote the i-th
 component of $d$. It is easy to check that $d_{0}$ and $d_{1}$
 coincide with the previously defined $d_{0}$ and $d_{1}$. Moreover
 $d_{i}(fg)=\sum_{j+k=i}d_{j}(f)d_{k}(g)$ for $f,g\in A$ since
 $d$ is an algebra homomorphism.
 \par We claim that $d_{i}$ restricted to $I^{i}$ is the natural
 homomorphism \\$I^{i}\to I^{i}/I^{i+1}$ for $i\geq {0}$. 
 This is true for $i=0$ and $i=1$. Let $f\in I^{i-1}$ and $g\in I$. By  
 induction
 \begin{equation*}
 \begin{aligned}
 d_{i}(fg)&=d_{i-1}(f)d_{1}(g)=(\text{ image of }f \text{ in }
 I^{i-1}/I^{i})(\text{ image of g }\text{ in }I/I^{2})\\ 
 &=\text{ image of }fg \text{ in }I^{i}/I^{i+1}. 
 \end{aligned}
 \end{equation*}
 This proves the claim since $I^{i}$ is additively generated by the elements
 of the form $fg$ with $f\in I^{i-1}\text{ and }g\in I$.
 \par We put $J=\bigoplus_{i\geq 1} S^{i}(I/I^{2})$ which is a $G$-invariant
 ideal of $S(I/I^{2})$.\linebreak Obviously for all 
 $n$\quad $d(I^{n})\subset J^{n}$.
 It follows that the induced homomorphism $A/I^{n}
 \to S(I/I^{2})/J^{n}$ is injective for all $n$.\par To prove its surjectivity
 it suffices to show that the homogenous elements of 
 $\bigoplus_{i\leq n-1}I^{i}/I^{i+1}$ differ from the image of $d$ by
 the elements of $\bigoplus_{i\geq n}I^{i}/I^{i+1}$.
 Let $x\in \bigoplus_{i\leq n-1}I^{i}/I^{i+1}$ with deg $x=k\quad (k\leq n-1)$.
 Then there exists $f\in I^{k}$ such that $d_{k}(f)=x$. We shall define
 inductively a sequence of elements $\{f_{i}\}_{k\leq i\leq n-1}$ with
 $f_{i}\in I^{i}$. Put $f_{k}=f$. Suppose the sequence
 $f_{k},f_{k+1},\dots f_{i}$ with $f_{j}\in I^{j}\quad k\leq j\leq i$
 has been defined. Let $f_{i+1}$ be an element of $I^{i+1}$ such that
 $d_{i+1}(f_{i+1})=-d_{i+1} (f_{k}+f_{k+1}+\dots + f_{i})$.
 Then \linebreak 
 $x-d(\sum_{i=k}^{n-1}f_{i})\in \bigoplus_{i\geq n}I^{i}/I^{i+1}$.
 So the induced map $A/I^{n}\to S(I/I^{2})/J^{n}$ is an isomorphism for  
 all $n$.
 \par We denote with $U$ the set of all the closed points of
 $Spec\, S(I/I^{2})$
 at which the morphism induced by $d \; Spec\, S(I/I^{2})\to Spec\, A$
 is etale.
 Then by  Proposition \ref{P:3}\quad $U$ is a $G$-invariant open subset of
 $E$ containing its zero section
 where $E$ denotes the normal bundle of $X$ in  $A^{n}_{k}$
 i.e the set of the closed points of $Spec\, S(I/I^{2})$.
 Moreover $Spec\, S(I/I^{2})\to Spec\, A$ maps $U$ into $A^{n}_{k}$
 which is the set of the closed points of $Spec\, A$. The restriction of the
 obtained $\phi\colon U\to A^{n}_{k}$ to the zero section of $E$ is an
 identity on $X$ 
 since
 $d\colon A \to S(I/I^{2})$ induces the identity homomorphism
 $R=A/I\to S(I/I^{2})/J=R$. Thus $\phi\colon U\to A^{n}_{k}$ has all the
 required properties.
 \end{proof}
 \begin{Rem}                                            
 The existence of $d\colon A\to S(I/I^{2})$
 which induces the isomorphism
 $A/I^{n} \simeq S(I/I^{2})/J^{n}$ for all $n$ has already
 been proved in the "absolute" case in \cite{B:2}.
 The proof actually shows that a $k$-algebra homomorphism
 \linebreak
 $d=(d_{i})\colon A \to S(I/I^{2})$ induces  the isomorphism
 $A/I^{n} \simeq S(I/I^{2})/J^{n}$ for all $n$ if
  $d_{0} \colon A \to A/I$
 is the natural homomorphism and the restriction of $d_{1} \colon A \to
 I/I^{2}$ to $I$ is the natural homomorphism $I\to I/I^{2}$. 
 \end{Rem}
 \begin{Rem}
 One easily obtains the following version of the Main Theorem in case
 the action of $G$ is not necessarily linear. \\
 Let $G$ be a reductive group acting (rationally) on an affine variety $Y$
 and let $X\subset Y$ be a closed, affine and smooth $G$-invariant subvariety. 
  Then there   exist
 (an algebraic) $G$-vector bundle $E$ on $X$, a $G$-invariant locally closed
 $U\subset E$ which
 contains the zero section of $E$ and $\phi\colon U\to Y$  such that
 $\phi $ is etale, $G$-invariant and  its restriction to the zero section
 of $E$ is the inclusion of $X$ in $Y$. \\Proof: $Y$ is isomorhic to a $G$-
 invariant closed affine subvariety of $A^{n}_{k}$ with a linear action
 of $G$ (\cite{K:1}). 
 Then by the Main Theorem there
 exist
 (an algebraic) $G$-vector bundle $E$ on $X$, a $G$-invariant open
 $U_{1}\subset E$ which
 contains the zero section of $E$ and $\phi\colon U_{1}\to A^{n}_{k}$ 
 such that
 $\phi $ is etale, $G$-invariant and  its restriction to the zero section
 of $E$ is the inclusion of $X$ in $A^{n}_{k}$. Put $U=\phi^{-1}(Y)$. Then
 $U\subset E$ and $\phi|U\colon U\to Y$ have all the required properties.

 \end{Rem}

\end{document}